\documentclass[a4paper]{article}

\usepackage[british]{babel}
\usepackage{amsmath}
\usepackage{amsfonts}
\usepackage{amsthm}
\usepackage{rotating}
\usepackage{verbatim}
\usepackage{enumerate}
\usepackage[cmtip,arrow]{xy}
\usepackage{pb-diagram,pb-xy}

\newtheorem{theorem}{Theorem}[section]

\newtheorem{proposition}[theorem]{Proposition}

\title{Transitive $A_6$-invariant $k$-arcs in $PG(2,q)$}
\author{ %
Massimo Giulietti
\\
\small{\texttt{giuliet@dmi.unipg.it}} \and G\'abor Korchm\'aros
\\
\small{\texttt{gabor.korchmaros@unibas.it}} \and
\\
Stefano Marcugini
\\
\small{\texttt{marcugini@dmi.unipg.it}}\and
\\
Fernanda Pambianco
\\
\small{\texttt{fernanda@dmi.unipg.it}}}
\date{}

\begin{document}
\maketitle

\begin{abstract}

For $q=p^r$ with a prime $p\ge 7$ such that $q \equiv 1$ or $19\pmod {30},$ the desarguesian
projective plane $PG(2,q)$ of order $q$ has a unique conjugacy class of projectivity groups
isomorphic to the alternating group $A_6$ of degree $6$. For a projectivity group
$\Gamma\cong A_6$ of $PG(2,q)$, we investigate the geometric properties of the (unique)
$\Gamma$-orbit $\mathcal{O}$ of size $90$ such that the $1$-point stabilizer of $\Gamma$ in
$\mathcal O$ is a cyclic group of order $4$. Here $\mathcal O$ lies either in $PG(2,q)$ or
in $PG(2,q^2)$ according as $3$ is a square or a non-square element in $GF(q)$. We show that
if $q\geq 349$ and $q\neq 421$, then $\mathcal O$ is a $90$-arc, which turns out to be
complete for $q=349, 409, 529, 601,661.$ Interestingly, $\mathcal O$ is the smallest known
complete arc in $PG(2,601)$ and in $PG(2,661).$ Computations are carried out by MAGMA.

\end{abstract}

\section{Introduction}

Let $GF(q)$ be a finite field of order $q=p^r$, a power of an odd prime $p$. In the
projective plane $PG(2,q)$ coordinatized by $GF(q)$, a $k-${\em arc} $K$ is a set of $k$
points no three of which are collinear. If an arc of $PG(2,q)$ is not contained in a larger
arc in $PG(2,q)$ then it is called {\em complete}. From the theory of linear codes, every
$k-$arc of $PG(2,q)$ corresponds to a $[k,3,k-2]$ {\em maximum distance separable} (MDS)
code of length $k$, dimension 3 and minimum distance $k-2$. This gives a motivation for the
the study of $k$-arcs in $PG(2,q)$; those with many projectivities were investigated in
several papers, see \cite{chk,CK,CK1,giu1,GG,hk1,hk2,hk3,kmp,KG1,kp,ks,P0}.

The maximum size of a (complete) arc in $PG(2,q)$ is $q+1$, and the points of an irreducible
conic in $PG(2,q)$ form an arc of size $q+1$. Actually, such $(q+1)$-arcs arising from
irreducible conics are the unique $(q+1)$-arcs in $PG(2,q)$. This is the famous Segre's
theorem \cite{segre1955}; see also \cite[Theorem 8.7]{H1}. Therefore, the projectivity group
which preserves a $(q+1)$-arc $K$ in $PG(2,q)$ is isomorphic to the projective linear group
$PGL(2,q)$ and acts on $K$ as $PGL(2,q)$ in its natural $3$-transitive permutation
representation. In particular, every $(q+1)$-arc $K$ is transitive. Here, the term of a
{\emph{transitive}} arc of $PG(2,q)$ is adopted to denote a $k$-arc $K$ such that the
projectivity group preserving $K$ acts transitively on the points of $K$.

Let $\Gamma$ be a finite group which can act faithfully as a projectivity group in $PG(2,q)$. Actually, this may happen in different characteristics $p.$ For instance,  $PG(2,q)$ with $p\neq 5$ has a projectivity group isomorphic to the alternating group $A_6$ if and only if $q\equiv 1$  or $19 \pmod {30}$, and in this case such a projectivity group is uniquely determined up to conjugacy in $PGL(3,q)$, see \cite{Bl}.  So  the question arises whether or not
a $\Gamma$-invariant arc of a fixed size $k$ exists in $PG(2,q)$ for infinitely many values of $p$.
{}From previous work, the answer is affirmative for $\Gamma\cong A_6$ and $k=72$, see \cite{kmp},
and $\Gamma\cong PSL(2,7)$ and $k=42$  see \cite{ik}. However the answer is negative for the
Hesse-group of order $216$ for any $k\geq 9$, see \cite{sonnino}.

In this paper we investigate the case of $\Gamma\cong A_6$ and $k=90$, giving a positive
answer to the above question:
\begin{theorem}
\label{mainth}
For a power $q$ of a prime $p\geq 7$ such that $q\equiv 1$ or $19 \pmod{30}$, let
$\Gamma\cong A_6$ be a projectivity group of $PG(2,q)$. Let $\mathcal{O}$ be the (unique)
$\Gamma$-orbit of length $90$ in $PG(2,q)$ such that
 the $1$-point stabilizer of $\Gamma$ in $\mathcal{O}$ is a cyclic group of order $4$.
Then $\mathcal{O}$ is a $90$-arc in $PG(2,q)$ except for a few cases where
\begin{itemize}

\item[\rm(i)] $q=61$ and $\mathcal{O}$ is a set of type $(0,1,2,4,6)$;
\item[\rm(ii)] $q=109$\,\,\, and $\mathcal{O}$ is a set of type $(0,1,2,3)$;
\item[\rm(iii)] $q=181$ and $\mathcal{O}$ is a set of type $(0,1,2,3)$;
\item[\rm(iv)] $q=229$\,\,\, and $\mathcal{O}$ is a set of type $(0,1,2,4)$;
\item[\rm(v)] $q=241$\,\,\, and $\mathcal{O}$ is a set of type $(0,1,2,4)$;
\item[\rm(vi)] $q=421$ and $\mathcal{O}$ is a2 set of type $(0,1,2,3)$;
\item[\rm(vii)] $q=7^2$ and $\mathcal{O}$ is a set of type $(0,1,2,4)$;
\item[\rm(viii)] $q=11^2$\,\,\, and $\mathcal{O}$ is a set of type $(0,1,2,5)$;
\item[\rm(ix)] $q=13^2$\,\,\, and $\mathcal{O}$ is a set of type $(0,1,2,4)$;
\item[\rm(x)] $q=17^2$\,\,\, and $\mathcal{O}$ is a set of type $(0,1,2,3)$;
\item[\rm(xi)] $q=19^2$\,\,\, and $\mathcal{O}$ is a set of type $(0,1,2,5)$;

\end{itemize}
\end{theorem}
An exhaustive computer aided search shows that such a $90$-arc may be complete for some
particular values of $q$, namely $q=349, 409, 529, 601, 661.$ It is worth mentioning that
this gives the smallest known complete arc in $PG(2,601)$ and in $PG(2,661),$ see
\cite{BDFMP, DGMP}.

Notation and terminology are standard, see \cite{H1}. Furthermore, $q$ always denotes a
power of an odd prime $p\geq 7$ such that $q\equiv 1$ or $19 \pmod{30}$.  Then $3$ divides $q-1$ and $5$ is a square element in
the multiplicative group of $GF(q)$. The latter two requirement are indeed necessary and
sufficient for $PGL(3,q)$ to have a subgroup $\Gamma \cong A_6$.

\section{Preliminary Results}
\label{sectgroup} We give an explicit representation of $\Gamma$ as a subgroup of $PGL(3,q)$ using the well known isomorphism $A_6\cong PSL(2,9)$.
Following \cite{kmp}, we choose a primitive element $\eta$ in $GF(9)$
satisfying $\eta^2 = \eta + 1$, and introduce the following matrices over $GF(9),$
\begin{displaymath}
\textsf{U}_1=\left(\begin{array}{cc}
1& 1 \\
0& 1
\end{array}\right),\,
\textsf{U}_2=\left(\begin{array}{cc}
 1& \eta^2   \\
0& 1
\end{array}\right),\,
\textsf{V}=\left(\begin{array}{cc}
 0 & -1  \\
1& 0
\end{array}\right),\,
\textsf{W}=\left(\begin{array}{cc}
 \eta & 0  \\
0& \eta^-1
\end{array}\right).
\end{displaymath}
It is easy to show that the above matrices generate $SL(2,9)$. Furthermore,
$\textsf{V}^4$ is the identity matrix $\textsf{I}$.

The factor group $SL(2,9)/\langle -\textsf{I}\rangle$ is $PSL(2,9).$

Let $\Phi:\, SL(2,9)\to PSL(2,9)$ be the associated natural homomorphism, and set
$M=\Phi(\textsf{M})$ with $\textsf{M}\in SL(2,9)$.\\
 There is a unique conjugacy class of
elements of order $4$ in $PSL(2,9)$, and the projectivity $W$ with matrix representation
$\textsf{W}$ is such an element of order $4$ (then $\langle W\rangle$ has order $4$ ...si
potrebbe aggiungere).

Now, fix a primitive third root $t$ of unity in $GF(q)$ and an element $z$ such that
$z^2=5$. Let $\Delta=t-t^2. $
Define the following matrices over $GF(q)$:\\
\begin{displaymath}
\mathbf{U}=\left(\begin{array}{lll}
0& 0 & 1  \\
1& 0 & 0 \\
0& 1 & 0
\end{array}\right),\,\,
\mathbf{\Omega}=\left(\begin{array}{lll}
 1&0 & 0  \\
 0& t & 0 \\
0& 0&  t^2
\end{array}\right),\,\,
\end{displaymath}
\vspace{0.2cm}
\begin{displaymath}
\mathbf{V}=\left(\begin{array}{lll}
-2& 1+ \Delta z & 1+\Delta z  \\
1-\Delta z & \;\;\;\; 4 & \; -2 \\
1- \Delta z &\; -2 &\;\;\; \; 4
\end{array}\right),\,\,
\mathbf{W}=\left(\begin{array}{lll}
 1& 1 & 1  \\
1& t  & t^2   \\
1& t^2   & t
\end{array}\right).
\end{displaymath}
Let $\bar{U},\bar{\Omega},\bar{V}$ and $\bar{W}$ be the associated projectivities of
$PGL(3,q)$. {}From \cite[Theorem 2.6]{kmp}, the projectivity  group generated by $\bar{U},
\bar{\Omega}, \bar{V}$ and $\bar{W}$ is isomorphic to $PSL(2,9)$. More precisely, the map
$\varphi$ with
\begin{equation}
\label{fai} \varphi:=\left\{
\begin{array}{lll}
U_1\to \bar{U}\\
U_2\to \bar{\Omega}\\
V\to \bar{V}\\
W\to \bar{W}
\end{array}
\right.
\end{equation} extends to an isomorphism from $PSL(2,9)$ into $PGL(3,q)$.
Therefore, the group generated by $\bar{U}, \bar{\Omega}, \bar{V}$ and $\bar{W}$ is
taken for $\Gamma$; that is,
\begin{equation}
\label{degGamma} \Gamma=\langle \bar{U}, \bar{\Omega}, \bar{V}, \bar{W}\rangle.
\end{equation}
A representative system of the $90$ cosets of $\langle W \rangle$ in $PSL(2,9)$ is listed
below.
\begin{equation}
\label{cos90}
\begin{array}{lll}
\{I, VWUW, WUV \Omega  VU, \Omega  VWU, UWUV  \Omega   U^2, UWV
\Omega U^2, V \Omega V \Omega ,
WV \Omega VU, \\
W^2 \Omega UV \Omega , VW^2 \Omega U, \Omega VUV \Omega , W^2V
\Omega UVU,
V \Omega VW^2UW,  \Omega VU^2VU,\\
WV \Omega ^2UV \Omega ,
 VU,  \Omega ^2VU^2VU, V \Omega VU^2VU,
V \Omega VW^2U, V \Omega U^2, U^2V \Omega ^2VU^2, \\
 U^2,  \Omega VW \Omega V \Omega , WVUV \Omega ,
V \Omega ^2VU^2, W \Omega V \Omega ,  \Omega ^2UVW^2U, WUWV \Omega UV, \\
 V \Omega ^2VUV, UW^2V \Omega ,
 \Omega VU, WUWUV \Omega , VU^2V \Omega UV,  \Omega ^2UVU,
 UWV \Omega ^2VU^2,\\
 W^2UWVUV,
VWV \Omega VU^2,  \Omega ^2VUV \Omega , WU^2V \Omega VU,  \Omega
^2U^2, VW^2 \Omega V \Omega U,
\\
UVU^2V \Omega ^2, UVW \Omega VU, W^2UV \Omega ^2V, VWU^2V \Omega ^2,
\Omega UVW^2U, W^2UWV^2,
V \Omega ,  \\
\Omega UVW \Omega , W^2UVU^2,W^2V^2U, V \Omega VW \Omega , V \Omega
^2U^2,VWV \Omega VUV,
 V \Omega UVW, \\
UW^2V^2 \Omega , VWV \Omega V, W^2VUV, UVW^2, UWV, UV, V, V \Omega
UVW^2, W^2V^2 \Omega , \\ \Omega ^2V \Omega , W^2V \Omega ^2UV, WV
\Omega VW, W^2V \Omega UV,
W^2VU^2, V \Omega VW^2, W^2V \Omega ^2U, \\
UWV^2U,  \Omega UVUWU^2,  \Omega ^2VU^2,
 \Omega VWV^2 \Omega U,  \Omega VW^2UW, V^2U^2V \Omega ^2,\\
 UV \Omega VW^2, \Omega VWUV \Omega ,
 \Omega VWV \Omega V, W \Omega UVU,  \Omega VW \Omega VUV, U^2VWU, \\
  W^2UV \Omega , VWV \Omega U,
   UV \Omega V,
 W^2 \Omega UV, V \Omega V \Omega U^2V, WUWUV,  \Omega UVW^2
\}.
\end{array}
\end{equation}
Replacing $U,V,W,\Omega$ with $\bar{U},\bar{V},\bar{W}.\bar{\Omega}$ gives a representative system of $\langle \bar{W}\rangle$ in $\Gamma$.

\section{The fixed points of $\bar{W}$}
The characteristic polynomial of $\bf{W}$ is $(\lambda^2-3)(\lambda-(1+2t))$ which has three
pairwise distinct roots, as $p\neq 3$. Let $s$ be an element in $GF(q)$ or in a quadratic
extension $GF(q^2)$ such that $s^2=3$. Then
$$
{\bf{v}}_1=\big(1,\frac{1}{2}(s-1),\frac{1}{2}(s-1)\big),\quad
{\bf{v}}_2=\big(1,-\frac{1}{2}(s+1),-\frac{1}{2}(s+1)\big),\quad {\bf{v}}_3=\big(0,1,-1\big)
$$
are three independent eigenvectors of $\bf{W}$. For $i=1,2,3$, let $P_i$ be the point represented by
${\bf{v}}_i$. Then $P_i$ are the fixed points of $\bar{W}$
in $PG(2,q)$ (or in $PG(2,q^2)$ when $s\in GF(q^2)\setminus GF(q)$).
The subgroup $S_2$ of $\Gamma$ generated by $\bar{V}$ and $\bar{W}$ is a dihedral group of order $8$.
Since $\bar{V}$ fixes $P_3$, this shows that $S_2$ is contained in the
stabilizer of $P_3$ in the action of $\Gamma$. But this is not consistent with the hypothesis
on the $1$-point stabilizer in Theorem \ref{mainth}. Furthermore, $\bar{V}$ interchanges the points
$P_1$ and $P_2$. Therefore, the $\Gamma$-orbit of $P_1$ contains $P_2$. From
the classification of subgroups of $A_6$, every proper subgroup of $\Gamma$ containing $\bar{W}$ also
contains $\bar{V}$. From this, the stabilizer
of $P_1$ under the action of $\Gamma$ is the group of order $4$ generated by $\bar{W}$. So, from
now on we may limit ourselves to consider the $\Gamma$-orbit $\mathcal O$  of $P_1$. We stress that
$\mathcal O$ is in $PG(2,q)$ (or in $PG(2,q^2)$ when $s\not\in GF(q)$).
The $90$ points in $\mathcal O$ can be computed as the images of the $P_1=(1,\frac{1}{2}(s-1),\frac{1}{2}(s-1))$ by the projectivities in the list
in (\ref{cos90}) after replacing $U,V,W,\Omega$ with $\bar{U},\bar{V},\bar{W},\bar{\Omega}$. These points are listed below.
\vspace{0.3cm}

 \begin{small}

\noindent $(2,-s-1,-s-1);((-12*s-12)*z*t+(-6*s-6)*z-6*s-18,(6*z+6*s)*t-6*z+6*s,(6*z-6*s)*t+12*z),\\
(((6*s+18)*z+18*s+54)*t+(12*s+36)*z-36*s+36,((6*s+18)*z+18*s-18)*t+(-6*s-18)*z+18*s+54,
((6*s-18)*z+18*s+18)*t+(-24*s-36)*z-36);\\
((2*s+6)*z*t+(s+3)*z+9*s+3,(2*s+6)*z*t+(s+3)*z-9*s-15,(2*s-6)*z*t+(s-3)*z+3*s+3);\\
(((s+3)*z+9*s+3)*t+(2*s+6)*z+6*s-6,((-5*s-3)*z+3*s+9)*t+(-4*s-6)*z+6,((-2*s-6)*z-12)*t+(-s-3)*z-3*s-9);\\
((-2*s*z-6)*t+2*s*z+6*s+12,(-2*s-6)*z*t+(-s-3)*z+3*s-3,(-2*s*z+6)*t-4*s*z+6*s+18);\\
((12*s+12)*z*t+(6*s+6)*z+6*s-18,(12*z-24*s-36)*t-12*z+12*s,(12*z+24*s+36)*t+24*z+36*s+36);\\
((-12*s*z+36)*t-24*s*z+36*s+108,(-12*s*z-36)*t+12*s*z+36*s+72,(-12*s-36)*z*t+(-6*s-18)*z+18*s-18);\\
(((3*s+15)*z+9*s+9)*t+(-3*s+3)*z+3*s+9,((-3*s-3)*z+3*s-9)*t+(3*s+3)*z-9*s-9,
((6*s+6)*z+6*s-18)*t+(3*s+3)*z-3*s-27);\\
((6*z-6*s)*t+12*z,(6*z+6*s)*t-6*z+6*s,(-12*s-12)*z*t+(-6*s-6)*z-6*s-18);\\
(((6*s+6)*z+18*s+18)*t+(-6*s-6)*z-6*s+18,((-6*s+6)*z-6*s-18)*t+(6*s+30)*z-18*s-18,((6*s+6)*z+(6*s+54))*t+(12*s+12)*z-12*s+36);\\
(((18*s+54)*z+18*s-90)*t+36*s+36,((-54*s-54)*z-18*s-90)*t-36*s*z-36*s-72,((-18*s-54)*z+18*s-90)*t+(-18*s-54)*z+90*s-18);\\
((36*s+108)*z*t+(18*s+54)*z-162*s-270,((-18*s+54)*z-54*s-54)*t+(18*s-54)*z-54*s-54,((-18*s-54)*z+162*s+54)*t+(-36*s-108)*z);\\
(((6*s+6)*z+18*s+18)*t+(12*s+12)*z+12*s-36,((-6*s-30)*z-18*s-18)*t+(6*s-6)*z-6*s-18,((6*s+6)*z-18*s-18)*t+(12*s+12)*z-24*s);\\
(((18*s+18)*z-30*s+18)*t+-24*s-72,(-12*s-36)*t+(18*s+18)*z+30*s-18,((-18*s-18)*z+6*s-18)*t+36*z-24*s-36);\\
(-4*z*t-2*z-2*s,-4*z*t-2*z-2*s,(-4*s-4)*z*t+(-2*s-2)*z-2*s-6);\\
(((6*s+6)*z-18*s-18)*t+(-6*s-6)*z-6*s-54,((-6*s-30)*z+18*s+18)*t+(-12*s-24)*z+12*s,((6*s+6)*z+18*s+18)*t+(-6*s-6)*z-6*s+18);\\
(((108*s+108)*z+36*s+180)*t+72*s*z+72*s+144,((36*s+108)*z-36*s+180)*t+(36*s+108)*z-180*s+36,((-36*s-108)*z-36*s+180)*t+-72*s-72);\\
((36*z+36*s)*t-36*z-72*s-108,(36*z-36*s)*t+72*z-108*s-108,(36*s+36)*z*t+(18*s+18)*z+18*s-54);\\
((2*z+2*s)*t-2*z+2*s,(-4*s-4)*z*t+(-2*s-2)*z-2*s-6,(2*z-2*s)*t+4*z);\\
(((-6*s-6)*z-18*s-18)*t+(-12*s-12)*z-12*s+36,((6*s-6)*z-6*s-18)*t+(12*s+24)*z+12*s,((-6*s-6)*z+6*s-18)*t+(6*s+6)*z-18*s-18);
(-s-1,2,-s-1);\\
(((-6*s+18)*z+18*s+18)*t+(-30*s-18)*z+18*s+54,((-6*s-18)*z-18*s-54)*t+(6*s+18)*z-18*s+18,((-6*s-18)*z-54*s-18)*t+(-12*s-36)*z-36*s+36);\\
((12*s+36)*z*t+(6*s+18)*z-54*s-90,((-6*s+18)*z+18*s+18)*t+(-12*s+36)*z,((-6*s-18)*z-54*s-18)*t+(6*s+18)*z-54*s-18);\\
((12*z-12*s)*t+24*z-36*s-36,(12*s+12)*z*t+(6*s+6)*z+6*s-18,(12*z+12*s)*t-12*z-24*s-36);\\
((-2*s-6)*z*t+(-s-3)*z+3*s-3,(-2*s*z-6*s-12)*t-4*s*z-6*s-18,(-2*s*z+6*s+12)*t+2*s*z-6);\\
(((-3*s-3)*z+9*s+9)*t+(3*s+3)*z+3*s+27,((-3*s-3)*z-9*s-9)*t+(3*s+3)*z+3*s-9,((3*s+15)*z-9*s-9)*t+(6*s+12)*z-6*s);\\
(((-18*s-54)*z+(18*s-90))*t+(-18*s-54)*z-18*s-126,((18*s+54)*z+18*s-90)*t+-72*s-72,((54*s+54)*z-18*s-90)*t+(18*s+54)*z+18*s-18);\\
(((-36*s-108)*z-36*s+180)*t+-72*s-72,((36*s+108)*z-36*s+180)*t+(36*s+108)*z-180*s+36,((108*s+108)*z+36*s+180)*t+72*s*z+72*s+144);\\
((-6*s+6)*z*t+(-3*s+3)*z+3*s+9,((-3*s-3)*z+(3*s+27))*t+(-6*s-6)*z,((-3*s-3)*z+(15*s+27))*t+(3*s+3)*z+15*s+27);\\
((-4*z-6*s-6)*t-2*z-2*s,(-4*z+(6*s+6))*t-2*z+4*s+6,(2*s+2)*z*t+(s+1)*z+s-3); \\
(((3*s+15)*z-9*s-9)*t+(6*s+12)*z-6*s,((6*s+6)*z+12*s)*t+(3*s+3)*z+9*s+9,((-3*s-3)*z-3*s-27)*t+(-6*s-6)*z+6*s-18);\\
(((36*s+108)*z+108*s+36)*t+216*s+72,(-72*s-72)*t+(36*s-108)*z-36*s-36,((-36*s-108)*z-108*s-180)*t+(-36*s-108)*z+108*s+180);\\
(((-s-1)*z-3*s-3)*t+(s+1)*z+s-3,((-s-1)*z+(3*s+3))*t+(s+1)*z+s+9,((s+5)*z-3*s-3)*t+(2*s+4)*z-2*s);\\
((-12*s*z-36*s-72)*t-24*s*z-36*s-108,(-12*s-36)*z*t+(-6*s-18)*z+18*s-18,(-12*s*z+(36*s+72))*t+12*s*z-36);\\
(((-18*s-54)*z+54*s+162)*t+(-36*s-108)*z+216,((90*s+54)*z+54*s+162)*t+(18*s-54)*z+54*s+54,((-18*s-54)*z-54*s-162)*t+(-36*s-108)*z+108*s-108);\\
((-144*s*z+216*s+648)*t-72*s*z+216*s+432,(72*s+216)*z*t+(36*s+108)*z-108*s+108,(-144*s*z-216*s-648)*t-72*s*z-216);\\
(((6*s+6)*z-18*s-18)*t+(12*s+12)*z-24*s,((12*s+24)*z+12*s)*t+(6*s+30)*z+18*s+18,((-12*s-12)*z-12*s+36)*t+(-6*s-6)*z+6*s+54);\\
(((6*s+18)*z-18*s+18)*t+(12*s+36)*z+72,((6*s+18)*z-18*s-54)*t+(-6*s-18)*z-54*s-18,((6*s-18)*z-18*s-18)*t+(30*s+18)*z-18*s-54);\\
((-s-1)*t,2,(s+1)*t+s+1);
((36*z-72*s-108)*t-36*z+36*s,(36*z+(72*s+108))*t+72*z+108*s+108,(36*s+36)*z*t+(18*s+18)*z+18*s-54);\\
((-12*s-12)*z*t+(-6*s-6)*z-6*s-54,((6*s+6)*z-30*s-54)*t+(-6*s-6)*z-30*s-54,((-6*s+6)*z-6*s-18)*t+(-12*s+12)*z);\\
((12*s*z+(36*s+72))*t+24*s*z+36*s+108,(12*s*z-36*s-72)*t-12*s*z+36,(12*s+36)*z*t+(6*s+18)*z-18*s+18);\\
(((18*s-18)*z+18*s+54)*t+(-18*s-90)*z+54*s+54,((-18*s-18)*z+(54*s+54))*t+(18*s+18)*z+18*s+162,((-18*s-18)*z-18*s+54)*t+(-36*s-36)*z-72*s);\\
((-12*s-36)*z*t+(-6*s-18)*z+54*s+90,((6*s+18)*z+(54*s+18))*t+(-6*s-18)*z+54*s+18,((6*s-18)*z-18*s-18)*t+(12*s-36)*z);\\
(((-3*s-3)*z-9*s-9)*t+(-6*s-6)*z-6*s+18,((-3*s-3)*z+(9*s+9))*t+(-6*s-6)*z+12*s,((3*s+15)*z+(9*s+9))*t+(-3*s+3)*z+3*s+9);\\
(-216*s,(-108*s-324)*t+-108*s-324,(108*s+324)*t);\\
((-4*s-4)*z*t+(-2*s-2)*z-2*s-6,(2*z-2*s)*t+4*z,(2*z+2*s)*t-2*z+2*s);\\
(((-s+3)*z+3*s+3)*t+(-5*s-3)*z+3*s+9,((-s-3)*z-9*s-3)*t+(-2*s-6)*z-6*s+6,((-s-3)*z-3*s-9)*t+(s+3)*z-3*s+3);\\
((6*s+6)*z*t+(3*s+3)*z+3*s+27,(-6*s+6)*z*t+(-3*s+3)*z+3*s+9,(6*s+6)*z*t+(3*s+3)*z-15*s-27);\\
(-108*s-108,-108*s-108,216);\\
((12*s+36)*z*t+(6*s+18)*z-54*s-90,((-6*s-18)*z+54*s+18)*t+(-12*s-36)*z,((-6*s+18)*z-18*s-18)*t+(6*s-18)*z-18*s-18);\\
((2*z-2*s)*t+4*z,(-4*s-4)*z*t+(-2*s-2)*z-2*s-6,(2*z+2*s)*t-2*z+2*s);\\
(((216*s+648)*z+(648*s+1944))*t+(-216*s-648)*z+648*s-648,((-1080*s-648)*z+(648*s+1944))*t+(-864*s-1296)*z+1296,((216*s+648)*z-648*s-1944)*t+(-216*s-648)*z-1944*s-
648);\\
(((-18*s-18)*z+6*s-90)*t+(-18*s-18)*z+30*s-18,(-36*z+(24*s+36))*t+(-18*s-54)*z+30*s+18,(-12*s-36)*t+(18*s+18)*z+30*s-18);
(-108*s-108,216*t,(108*s+108)*t+108*s+108);\\
((72*s+216)*z*t+(36*s+108)*z-108*s+108,(-144*s*z-216*s-648)*t-72*s*z-216,(-144*s*z+(216*s+648))*t-72*s*z+216*s+432);\\
((-36*s-36)*z*t+(-18*s-18)*z+90*s+162,(36*s-36)*z*t+(18*s-18)*z-18*s-54,(-36*s-36)*z*t+(-18*s-18)*z-18*s-162);
((-6*s+6)*z*t+(-3*s+3)*z+3*s+9,(6*s+6)*z*t+(3*s+3)*z+3*s+27,(6*s+6)*z*t+(3*s+3)*z-15*s-27);\\
((-2*s-6)*z*t+(-s-3)*z+3*s-3,(4*s*z-6*s-18)*t+2*s*z-6*s-12,(4*s*z+(6*s+18))*t+2*s*z+6);\\
((-2*s+2)*z*t+(-s+1)*z+s+3,(2*s+2)*z*t+(s+1)*z-5*s-9,(2*s+2)*z*t+(s+1)*z+s+9);\\
((-4*s-4)*z*t+(-2*s-2)*z-2*s-6,-4*z*t-2*z-2*s,-4*z*t-2*z-2*s);\\
(((-18*s-54)*z+(18*s-90))*t+(-18*s-54)*z+90*s-18,((-54*s-54)*z-18*s-90)*t-36*s*z-36*s-72,((18*s+54)*z+(18*s-90))*t+36*s+36);
(216,(-108*s-108)*t,(108*s+108)*t+108*s+108);\\
((2*s+2)*z*t+(s+1)*z+s-3,(2*z-2*s)*t+4*z-6*s-6,(2*z+2*s)*t-2*z-4*s-6);\\
(((-18*s-54)*z-18*s+90)*t+72*s+72,((18*s+54)*z-18*s+90)*t+(18*s+54)*z+18*s+126,
((-54*s-54)*z+(18*s+90))*t+(-18*s-54)*z-18*s+18);\\
((-36*s-36)*z*t+(-18*s-18)*z+90*s+162,(-36*s-36)*z*t+(-18*s-18)*z-18*s-162,(36*s-36)*z*t+(18*s-18)*z-18*s-54);\\
(((-18*s-54)*z+(18*s-90))*t+(-18*s-54)*z+90*s-18,((18*s+54)*z+(18*s-90))*t+36*s+36,((-54*s-54)*z-18*s-90)*t-36*s*z-36*s-72);\\
(-12*z*t-6*z-6*s,(-12*s-12)*z*t+(-6*s-6)*z-6*s-18,-12*z*t-6*z-6*s);\\
((36*s+36)*z*t+(18*s+18)*z+18*s-54,(36*z+36*s)*t-36*z-72*s-108,(36*z-36*s)*t+72*z-108*s-108);\\
((6*z+6*s)*t-6*z+6*s,(6*z-6*s)*t+12*z,(-12*s-12)*z*t+(-6*s-6)*z-6*s-18);\\
((-36*s-108)*t+-36*s-108,(36*s+108)*t,-72*s);\\
(((-s-3)*z+3*s-3)*t+(-2*s-6)*z-12,((-s+3)*z+(3*s+3))*t+(-5*s-3)*z+3*s+9,((-s-3)*z+(3*s+9))*t+(s+3)*z+9*s+3);\\
((-4*z-6*s-6)*t-2*z-2*s,(2*s+2)*z*t+(s+1)*z+s-3,(-4*z+(6*s+6))*t-2*z+4*s+6);\\
(((-36*s-108)*z+(324*s+108))*t+(-72*s-216)*z,((-36*s-108)*z+(324*s+540))*t+(36*s+108)*z+324*s+
540,(72*s-216)*z*t+(36*s-108)*z+108*s+108);\\
((6*s-18)*z*t+(3*s-9)*z+9*s+9,((-3*s-9)*z+(27*s+45))*t+(3*s+9)*z+27*s+45,((-3*s-9)*z+(27*s+9))*t+(-6*s-18)*z);\\
((-72*s+72)*z*t+(-36*s+36)*z+36*s+108,((-36*s-36)*z-36*s-324)*t+(36*s+36)*z-36*s-324,((-36*s-36)*z-180*s-324)*t+(-72*s-72)*z);\\
(((18*s-18)*z-18*s-54)*t+(36*s+72)*z+36*s,((-18*s-18)*z+(18*s-54))*t+(18*s+18)*z-54*s-54,((-18*s-18)*z-54*s-54)*t+(-36*s-36)*z-36*s+108);\\
((-12*s-36)*z*t+(-6*s-18)*z-54*s-18,((6*s-18)*z-18*s-18)*t+(12*s-36)*z,((6*s+18)*z-54*s-90)*t+(-6*s-18)*z-54*s-90);\\
(((-36*s+108)*z+(108*s+108))*t+(-180*s-108)*z+108*s+
324,((72*s+216)*z+(216*s-216))*t+(36*s+108)*z-108*s-324,((72*s+216)*z+432)*t+(36*s+108)*z+108*s+324);\\
(((s+3)*z-3*s-9)*t+(2*s+6)*z-12,((s+3)*z+(3*s+9))*t+(2*s+6)*z-6*s+6,((-5*s-3)*z-3*s-9)*t+(-s+3)*z-3*s-3);\\
(((108*s+108)*z-180*s+108)*t+72*s+216,(144*s+432)*t+(108*s+108)*z+180*s-108,(216*z+(144*s+216))*t+(-108*s-108)*z-36*s+108);\\
((-4*s*z+(6*s+18))*t-2*s*z+6*s+12,(-4*s*z-6*s-18)*t-2*s*z-6,(2*s+6)*z*t+(s+3)*z-3*s+3);\\
((-6*s+6)*z*t+(-3*s+3)*z+3*s+9,((-3*s-3)*z-15*s-27)*t+(-6*s-6)*z,((-3*s-3)*z-3*s-27)*t+(3*s+3)*z-3*s-27);\\
((-36*s-108)*t,(36*s+108)*t+36*s+108,72*s);\\
(((6*s-6)*z-6*s-18)*t+(12*s+24)*z+12*s,((-6*s-6)*z-18*s-18)*t+(-12*s-12)*z-12*s+36,((-6*s-6)*z+(6*s-18))*t+(6*s+6)*z-18*s-18);\\
(((3*s+15)*z+9*s+9)*t+(-3*s+3)*z+3*s+9,((-3*s-3)*z+(9*s+9))*t+(-6*s-6)*z+12*s,((-3*s-3)*z-9*s-9)*t+(-6*s-6)*z-6*s+18);\\
((72*s*z-72*s-288)*t+72*s*z-144*s-360,(72*s-72)*t+(36*s+108)*z+36*s-36,(-72*s*z-72*s-288)*t+72*s+72);\\
(((3*s+15)*z-9*s-9)*t+(6*s+12)*z-6*s,((-3*s-3)*z-9*s-9)*t+(3*s+3)*z+3*s-9,((-3*s-3)*z+(9*s+9))*t+(3*s+3)*z+3*s+27);\\
(((3*s+15)*z+9*s+9)*t+(-3*s+3)*z+3*s+9,((-3*s-3)*z-9*s-9)*t+(-6*s-6)*z-6*s+18,((-3*s-3)*z+(9*s+9))*t+(-6*s-6)*z+12*s)).$

\vspace{0.3cm}
 \end{small}

Let ${\mathcal {O}}=\{P_1,Q_1,\ldots,Q_{89}\}$. The points $P_1,Q_i$ and $Q_j$ are collinear
if and only if the determinant $D_{i,j}$ of the coordinates of these points vanishes. There
are $3916$ triples $\{P_1,Q_i,Q_j\}$ with $1\le i<j\leq 89$. Observe that $D_{i,j}$ can be
viewed as a polynomial in $t,s$ and $z$, say $D_{i,j}(t,s,z)$, with coefficients in
$\mathbb{Z}$.  Therefore a necessary and sufficient  condition for the points $Q_i,Q_j\in
\mathcal {O}$ to produce together with $P_1$ a collinear triple  is that $(t,s,z)$ be a
solution of the system of equations
\begin{equation}
\label{system}\left\{
\begin{array}{lll}
t^2+t+1=0;\\
s^2=3;\\
z^2=5;\\
D_{i,j}(t,s,z)=0.
\end{array}
\right.
\end{equation}
We look at the above system over $\mathbb{Z}$ with unknowns $t,s,z$ and use Sylvester's
resultant to discuss solvability. Eliminating $t$ from the first and the forth equations
produces an equation in $s,z$ over $\mathbb{Z}$: then eliminating $s$ from this and the
second equation provides an equation in $z$ over $\mathbb{Z}$; finally eliminating $z$ from
this and the third equation gives an integer, the resultant of the system. A sufficient
condition for a triple of points not to be collinear is that this resultant does not vanish
in $Z_p$. \\
A computer aided search shows that such a resultant is a non zero integer for any of the
above $3916$ cases. Now, let $\delta_{i,j}$ be the set of all prime divisors of the
resultant arising from the triple $\{P_1,Q_i,Q_j\}$. If $p\not\in \delta_{i,j}$ then the
points $P_1,Q_i,Q_j$ are not collinear. More generally, let $\delta$ be the set of all
primes appearing in some of the $3916$ sets $\delta_{i,j}$. If $p\not\in \delta$ then the
$\Gamma$-orbit $\mathcal O$ is $90$-arc.


An exhaustive computer-aided computation shows that $\delta$ has size $14$, namely
$\delta=\{2,3,5,7,11,13,17,19,61,109,181,229,241,421\}.$ Therefore, the following result
holds.
\begin{proposition}
\label{th9giu} The $\Gamma$-orbit $\mathcal{O}$ of the point $P_1$ has length $90$ and the stabilizer of $P_1$ in $\Gamma$ is a cyclic group of order $4$. Furthermore, $\mathcal{O}$ is a $90$-arc on $PG(2,q)$ with
$q=p^h$ and $p\geq 7$ apart from finitely many values of $p$ which are
$$7,11,13,17,19,61,109,181,229,241,421.$$
\end{proposition} Now, we discuss the exceptional cases.
\subsection{$p=7, 11, 13, 17$ }

In this case $p^2\equiv 1$ or $19\pmod {30}$. Therefore, $\mathcal O$ lies in
$PG(2,p^2)\setminus PG(2,p)$. By a MAGMA \cite{mag} computation, some $\delta_{i,j}$ is
divisible by $p$. Hence $\mathcal{O}$ is not an arc. Some more effort allows to compute the
intersection numbers of $\mathcal{O}$ with lines. The results are reported below.
\begin{itemize}


\item For $q = 7^{2}$ a square root of 5 is $w^{20}$, where  $w$ is a
primitive element of $GF(7^2)$ such that $w^{2}+6w+3=0$. In this case $\mathcal{O}$ is a
complete $(90,4)-$arc with 336 external lines, 810 tangents, 765 bi-secants, 540
four-secants.

\item For $q = 11^{2}$ a primitive cubic root of unity is $w^{40}$, where
$w$ is a primitive element
 of $GF(11^2)$.
In this case $\mathcal{O}$ is a non-complete $(90,5)-$arc with 7248 external lines, 4320
tangents, 3105 bi-secants, 90 five-secants.

\item For $q = 13^{2}$ a square root of 5 is $w^{63}$, where  $w$ is a
primitive element of $GF(13^2)$
 such that $w^{2}+12w+2=0$.
In this case $\mathcal{O}$ is a non-complete $(90,4)-$arc with 16896 external lines, 8730
tangents, 2925 bi-secants, 180 four-secants.

\item For $q = 17^{2}$ a primitive cubic root of unity is $w^{98}$ and a
square root of 5 is $w^{45}$, where  $w$ is a primitive element of
$GF(17^2)$
 such that $w^{2}+16w+3=0$.
In this case $\mathcal{O}$ is a non-complete $(90,3)-$arc with 61356 external lines, 19170
tangents, 2925 bi-secants, 360 three-secants.

\end{itemize}

\subsection{$p=19$}

$\Gamma$ is a projectivity group of $PG(2,19)$. However,  $s\in GF(19^2)\setminus GF(19)$,
whence $\mathcal{O}$ lies in $PG(2,19^2)\setminus PG(2,19)$. Furthermore, $\mathcal O$ is a
non-complete $(90,5)-$arc with 101676 external lines, 25650 tangents, 3285 bi-secants, 72
five-secants. Here, $s=w^{130}$ where $ w^2+18w +2 =0$.

\subsection{$p=61, 109, 181, 229, 241, 421$ }
In this case, $\Gamma$ is a projective group of $PG(2,p)$ and $s\in GF(p)$. Therefore
$\mathcal{O}$ lies in $PG(2,p)$. Again, some $\delta_{i,j}$ is divisible by $p$, and
$\mathcal{O}$ is not an arc. By a MAGMA computation, the intersection numbers of
$\mathcal{O}$ lines can be computed, and the results are reported below.
\begin{itemize}

\item For $p=61$  $\mathcal{O}$ is a non-complete $(90,6)-$arc with 1068
external lines, 450 tangents, 2025 bi-secants, 180 four-secants, 60 six-secants.

\item For $p=109$ $\mathcal{O}$ is a non-complete $(90,3)-$arc
with 5736 external lines, 2970 tangents, 2925 bi-secants, 360 three-secants.

\item For $p=181$  $\mathcal{O}$ is a non-complete $(90,3)-$arc
with 20208 external lines, 9450 tangents, 2925 bi-secants, 360 three-secants.

\item For $p=229$  $\mathcal{O}$ is a non-complete $(90,4)-$arc
with 35436 external lines, 14130 tangents, 2925 bi-secants, 180 four-secants.

\item For $p=241$ $\mathcal{O}$ is a non-complete $(90,4)-$arc
with 40008 external lines, 15210 tangents, 2925 bi-secants, 180 four-secants.

\item For $p=421$  $\mathcal{O}$ is a non-complete $(90,3)-$arc with
143328 external lines, 31050 tangents, 2925 bi-secants, 360 three-secants.

\end{itemize}

The results of the present section provide a proof of Theorem 1.1.

\end{document}